\documentclass[11pt,a4paper]{article}
\usepackage[english]{babel}
\usepackage{amsmath,amsfonts,amsthm,amssymb,amsbsy,upref,color,graphicx,amscd,hyperref,makeidx,enumerate,bbm,comment}
\usepackage[a4paper,margin=3.0truecm]{geometry}
\usepackage[active]{srcltx}
\usepackage[latin1]{inputenc}

\def\ds{\displaystyle}
\def\eps{{\varepsilon}}

\def\O{\Omega}
\def\R{\mathbb{R}}

\def\A{\mathcal{A}}
\def\E{\mathcal{E}}

\def\HH{\mathcal{H}}

\newcommand{\be}{\begin{equation}}
\newcommand{\ee}{\end{equation}}
\newcommand{\bib}[4]{\bibitem{#1}{\sc#2: }{\it#3. }{#4.}}

\numberwithin{equation}{section}
\theoremstyle{plain}
\newtheorem{teo}{Theorem}[section]

\theoremstyle{remark}

\newenvironment{ack}{{\bf Acknowledgements. }}

\title{Two optimization problems in thermal insulation}

\author{Dorin Bucur, Giuseppe Buttazzo, Carlo Nitsch}

\begin{document}

\maketitle

\begin{abstract}
We consider two optimization problems in thermal insulation: in both cases the goal is to find a thin layer around the boundary of the thermal body which gives the best insulation. The total mass of the insulating material is prescribed.. The first problem deals with the case in which a given heat source is present, while in the second one there are no heat sources and the goal is to have the slowest decay of the temperature. In both cases an optimal distribution of the insulator around the thermal body exists; when the body has a circular symmetry, in the first case a constant heat source gives a constant thickness as the optimal solution, while surprisingly this is not the case in the second problem, where the circular symmetry of the optimal insulating layer depends on the total quantity of insulator at our disposal. A symmetry breaking occurs when this total quantity is below a certain threshold. Some numerical computations are also provided, together with a list of open questions.
\end{abstract}



\section{Introduction}\label{sintro}

Insulation problems will represent some of the crucial research fields in the next future, being related to energy saving, pollution control, environment improving. Several sciences are involved in this topic: Civil Engineering for the design of new buildings with more efficient energy consumption, Physics and Chemistry for the research on new material with better insulating properties, Mathematics for the study of the partial differential equations involved in the heat conduction in the presence of insulating regions.

In this short note we present two problems related to optimal insulation; in both cases a domain $\O$ of $\R^d$ is given, we assume it is a bounded open set with a regular boundary $\partial\O$. For simplicity we assume that $\O$ is a conducting domain with a constant conductivity coefficient, that we assume equal to one. The goal is to distribute around $\partial\O$ a layer $\Sigma$ of insulating material in some efficient way; the optimality criterion we use is described later. We describe the layer $\Sigma$ by means of the tangential and normal coordinates on $\partial\O$:
$$\Sigma_\eps=\big\{\sigma+t\nu(\sigma)\ :
\ \sigma\in\partial\O,\ 0\le t<\eps h(\sigma)\big\},$$
where $\nu(\sigma)$ is the exterior normal versor to $\partial\O$ at the point $\sigma$ and the function $h$ describes the variable thickness. The index $\eps$ describes the average thickness of the layer and is taken very small (for a house of several meters of diameter the thickness of the insulating layer is usually of few centimeters). The conductivity coefficient of the insulating material in the layer $\Sigma_\eps$ is taken very small too, we denote it by $\delta$. Finally, we assume that the temperature is zero outside the set $\O\cup\Sigma_\eps$

The two problems we deal with are described below in a precise mathematical form.

\bigskip

{\bf Problem 1. }We put in $\O$ a heat source $f\in L^2(\O)$; after waiting enough time the temperature $u(t,x)$ approaches the stationary temperature $u(x)$ that is the solution of the elliptic equation with transmission conditions at $\partial\O$
\be\label{pde1}
\begin{cases}
-\Delta u=f&\hbox{in }\O\\
-\Delta u=0&\hbox{in }\Sigma_\eps\\
u=0&\hbox{on }\partial(\O\cup\Sigma_\eps)\\
\ds\frac{\partial u^-}{\partial\nu}
=\delta\frac{\partial u^+}{\partial\nu}&\hbox{on }\partial\O\;.
\end{cases}
\ee
Equivalently, the stationary temperature $u$ can be seen as the solution of the minimum problem on $H^1_0(\O\cup\Sigma_\eps)$ for the energy functional
\be\label{cost1}
E_{\eps,\delta}(u)=\frac12\int_\O|\nabla u|^2\,dx+\frac\delta2\int_{\Sigma_\eps}|\nabla u|^2\,dx-\int_\O fu\,dx\;.
\ee
The PDE \eqref{pde1} is indeed the Euler-Lagrange equation of the minimum problem related to the cost functional \eqref{cost1}. Denoting by $u$ the solution of the PDE \eqref{pde1} (or of the variational problem for the energy \eqref{cost1}) a multiplication by $u$ in \eqref{pde1} and a standard integration by parts allows us to write the minimum of the energy functional \eqref{cost1} in the form
$$\min_{H^1_0(\O\cup\Sigma_\eps)}E_{\eps,\delta}=-\frac12\int_\O fu\,dx\;.$$
Note that, when the heat sources are uniformly distributed, the minimization of the energy functional above corresponds to the maximization of the average temperature.

\begin{figure}[htp]
\begin{center}
\includegraphics[height=6truecm]{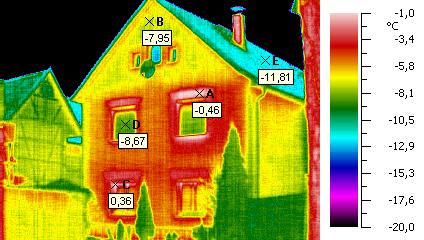}
\caption{The thermography of a house (photo by Lutz Weidner on commons.wikimedia.org).}
\label{fig1a}
\end{center}
\end{figure}

\begin{figure}[htp]
\begin{center}
\includegraphics[height=6truecm]{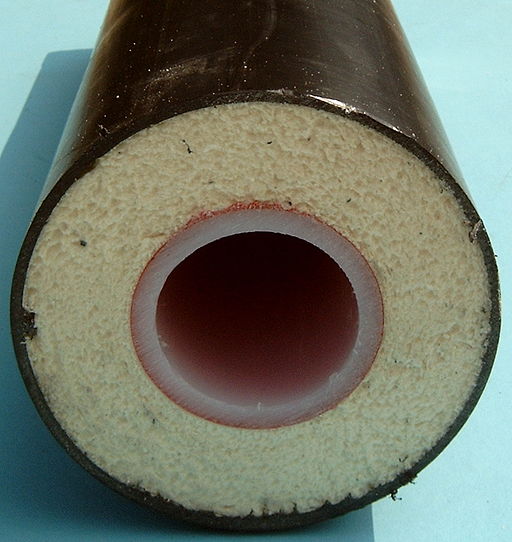}
\caption{The thermal insulation of a pipe (photo by S\"onke Kraft aka Arnulf zu Linden on commons.wikimedia.org).}
\label{fig1b}
\end{center}
\end{figure}

\begin{figure}[htp]
\begin{center}
\includegraphics[height=6truecm]{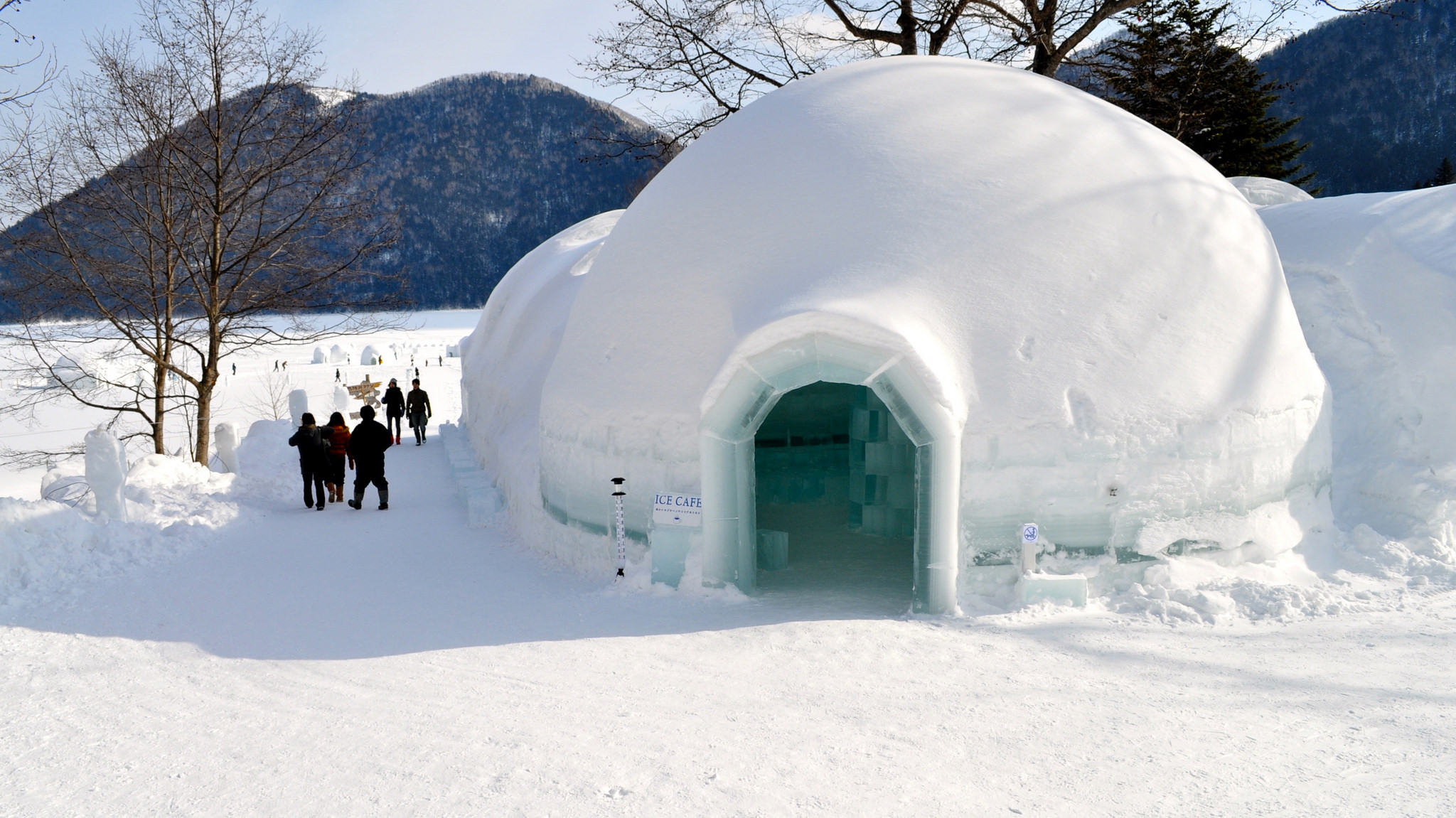}
\caption{The thermal insulation of an igloo by a layer of snow (photo by David McKelvey on www.flickr.com).}
\label{figigloo}
\end{center}
\end{figure}

The optimization problem we deal with consists in the optimal choice of the shape of the insulating layer $\Sigma_\eps$ around $\partial\O$ once the total amount of insulating material is prescribed. Stressing the dependence on $h$ of the energy functional and denoting by $E(h)$ the quantity
$$E(h)=\min_{H^1_0(\O\cup\Sigma_\eps)}E_{\eps,\delta}$$
our first optimization problem can be written in the form
$$\min\big\{E(h)\ :\ h\in\HH_m\big\}$$
where $\HH_m$ denotes the class of admissible choices
\be\label{hm}
\HH_m=\left\{h:\partial\O\to\R,\ h\ge0,\ \int_{\partial\O}h\,d\HH^{d-1}=m\right\}.
\ee
This applies for instance to the thermal insulation of an house (see Figure \ref{fig1a}) or of a pipe (see Figure \ref{fig1b}). One of the most crucial questions in this field is: ``which parts have to be more protected?'' For instance, in a radially symmetric body (like an igloo, see Figure \ref{figigloo}), should we put a layer of insulating material with a constant boundary thickness?

\begin{figure}[htp]
\begin{center}
\includegraphics[height=6truecm]{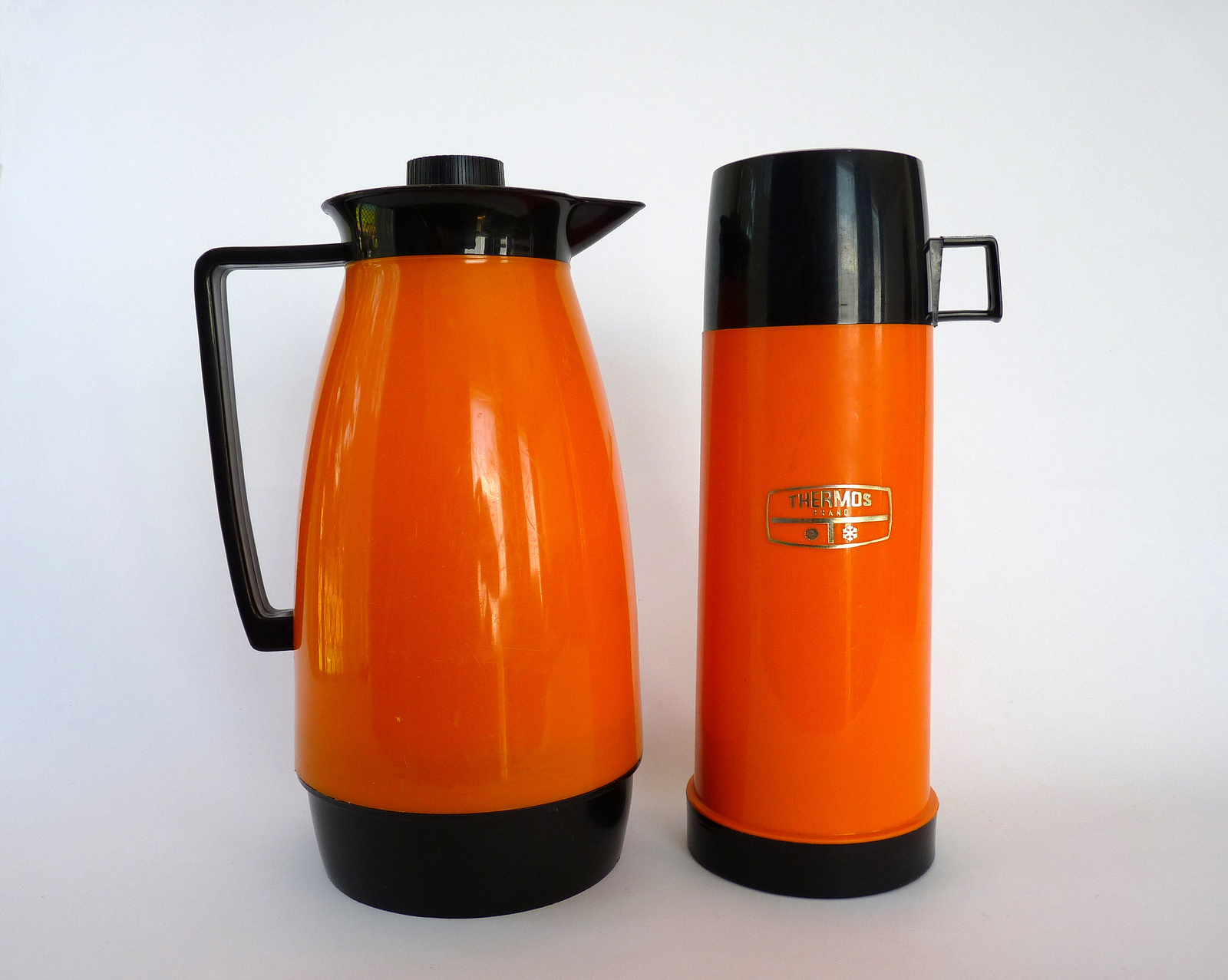}
\caption{The thermal insulation of a coffee pot (photo by Gail Thomas on www.flickr.com).}
\label{fig2}
\end{center}
\end{figure}

{\bf Problem 2. }The second problem that we consider deals with a domain $\O$ as above, with a fixed initial temperature $u_0$ and without any heat source. In this case the temperature decays to zero and our goal is to put the insulating material around $\O$ in order this decay be as low as possible. This applies for instance to the thermal insulation of a coffee pot, see Figure \ref{fig2}.

By the Fourier analysis of the corresponding heat diffusion equation, the decay of the temperature goes as $e^{-t\lambda}$, where $\lambda$ is the first eigenvalue of the elliptic operator written in the weak form as
$$\langle\A u,\phi\rangle=\int_\O\nabla u\nabla\phi\,dx+\delta\int_{\Sigma_\eps}\nabla u\nabla\phi\,dx\qquad\forall\phi\in H^1_0(\O\cup\Sigma_\eps)\;.$$
Therefore, indicating by $\lambda(h)$ the first eigenvalue above, stressing its dependence on the function $h$, Problem 2 reduces to the minimization problem
\be\label{pb2}
\min\big\{\lambda(h)\ :\ h\in\HH_m\big\}\;,
\ee
where $\HH_m$ is the class of admissible choices introduced in \eqref{hm}.

\section{The asymptotic problem}\label{sasympt}

In order to simplify the questions illustrated in Section \ref{sintro} we consider the asymptotic model when both the thickness $\eps$ of the insulating layer, as well as the conductivity coefficient $\delta$ of the insulator tend to zero. The identification of the limit problem goes back to \cite{sp74} (see also \cite{sp80}), a more general variational framework was considered in \cite{brca80} and its formulation in terms of $\Gamma$-convergence can be found in \cite{acbu86}. We summarize the results in this last point of view by considering the functionals
$$\frac12\int_\O|\nabla u|^2\,dx+\frac\delta2\int_{\Sigma_\eps}|\nabla u|^2\,dx\qquad u\in H^1_0(\O\cup\Sigma_\eps)\;.$$

\begin{itemize}
\item When $\eps\ll\delta$ the limit problem is the Dirichlet one, related to the functional
$$\frac12\int_\O|\nabla u|^2\,dx\qquad u\in H^1_0(\O)\;.$$

\item When $\eps\gg\delta$ the limit problem is the Neumann one, related to the functional
$$\frac12\int_\O|\nabla u|^2\,dx\qquad u\in H^1(\O)\;.$$

\item When $\eps\approx k\delta$ with $k>0$ the limit problem is a Robin type problem, related to the functional
$$\frac12\int_\O|\nabla u|^2\,dx+\frac{1}{2k}\int_{\partial\O}\frac{u^2}{h}\,d\sigma\qquad u\in H^1(\O)\;.$$
\end{itemize}

In the rest of this article we are in the framework of the last situation. We can reformulate now Problems 1 and 2 in their asymptotic form, as $\eps$ and $\delta$ go to zero, with $k\approx\eps/\delta$.

\medskip

The asymptotic form of Problem 1 becomes now
\be\label{pb1}
\min\big\{\E(h)\ :\ h\in\HH_m\big\}
\ee
where $\HH_m$ is given in \eqref{hm} and $\E$ is the asymptotic energy
$$\E(h)=\min\left\{\frac12\int_\O|\nabla u|^2\,dx+\frac{1}{2k}\int_{\partial\O}\frac{u^2}{h}\,d\sigma-\int_\O fu\,dx\ :\ u\in H^1(\O)\right\}\;.$$
The Euler-Lagrange equation of the minimum problem above is
$$\begin{cases}
-\Delta u=f\quad\hbox{in }\O\\
\ds\frac1k u+h\frac{\partial u}{\partial\nu}=0\quad\hbox{on }\partial\O\;.
\end{cases}$$
Denoting by $u_h$ its solution, multiplicating both sides by $u_h$ and integrating by parts gives that
$$\E(h)=-\frac12\int_\O fu_h\,dx\;.$$

Analogously, the asymptotic form of Problem 2 is given by
\be\label{pb2}
\min\big\{\lambda(h)\ :\ h\in\HH_m\big\}
\ee
where $\lambda(h)$ is the first eigenvalue of the elliptic operator written in a weak form as
$$\langle\A u,\phi\rangle=\int_\O\nabla u\nabla\phi\,dx+\frac1k\int_{\partial\O}\frac{u\phi}{h}\,d\sigma\qquad\forall\phi\in H^1(\O)\;.$$
Equivalently, $\lambda(h)$ can be written in terms of the Rayleigh quotient
$$\lambda(h)=\min\left\{\frac{\ds\int_\O|\nabla u|^2\,dx+\frac1k\int_{\partial\O}\frac{u^2}{h}\,d\sigma}{\ds\int_\O u^2\,dx}\ :\ u\in H^1(\O),\ u\ne0\right\}\;.$$

\section{Energy optimization}\label{senergy}

The optimization problem \eqref{pb1} is a double minimization problem:
$$\min_{h\in\HH_m}\,\min_{u\in H^1(\O)}\left\{\frac12\int_\O|\nabla u|^2\,dx+\frac{1}{2k}\int_{\partial\O}\frac{u^2}{h}\,d\sigma-\int_\O fu\,dx\right\};$$
interchanging the two minima we have that for every $u\in H^1(\O)$ which does not identically vanish on $\partial\O$ the best choice for $h$ is
$$h=m\frac{|u|}{\int_{\partial\O}|u|\,d\sigma}\;,$$
while the choice of $h$ is irrelevant when $u\in H^1_0(\O)$. This reduces the minimization problem \eqref{pb1} to
\be\label{pb1aux}
\min\left\{\frac12\int_\O|\nabla u|^2\,dx+\frac{1}{2km}\Big(\int_{\partial\O}|u|\,d\HH^{d-1}\Big)^2-\int_\O fu\,dx\ :\ u\in H^1(\O)\right\}.
\ee
The problem has been studied in \cite{bu88, bubuni}; we summarize here the results.

\begin{teo}\label{teoenergy}
Assume $\O$ is connected. Then the functional
$$F(u)=\frac12\int_\O|\nabla u|^2\,dx+\frac{1}{2km}\Big(\int_{\partial\O}|u|\,d\HH^{d-1}\Big)^2$$
is strictly convex on $H^1(\O)$, hence for every $f\in L^2(\O)$ the minimization problem \eqref{pb1aux} admits a unique solution $\bar u$. Thus the optimal function $h_m$ for problem \eqref{pb1} is given by
$$h=m\frac{|\bar u|}{\int_{\partial\O}|\bar u|\,d\sigma}\;.$$
\end{teo}

By uniqueness, if $\O=B_R$ in $\R^d$ and $f=1$ the optimal solution $\bar u$ above is radial:
$$\bar u(r)=\frac{R^2-r^2}{2d}+\frac{km}{d^2\omega_d R^{d-2}}\;,$$
hence, the optimal thickness $h_m$ is constant.

\begin{figure}[htp]
\begin{center}
\includegraphics[height=6truecm]{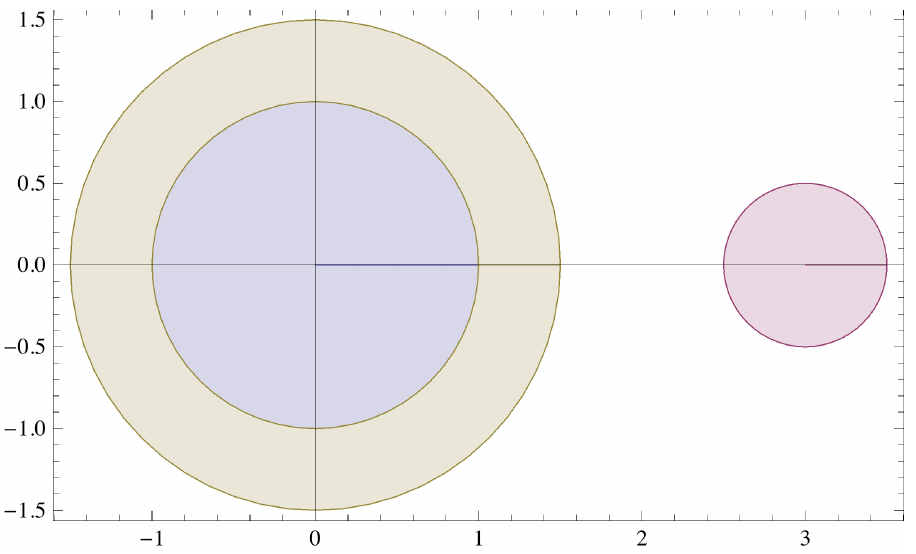}
\caption{Optimal insulation of a domain composed by two disjoint balls of different radius.}
\label{fig3}
\end{center}
\end{figure}

If $\O$ is not connected the optimal insulation strategy is different. Let $\O=B_{R_1}\cup B_{R_2}$ in $\R^d$ (union of two disjoint balls), and $f=1$. Then:
\begin{itemize}
\item if $R_1=R_2=R$ any choice of $h_m$ constant around $B_{R_1}$ and on $B_{R_2}$ is optimal;
\item if $R_1\ne R_2$ then the optimal choice is to concentrate all the insulator around the largest ball, with constant thickness, leaving the smallest ball unprotected (see Figure \ref{fig3}).
\end{itemize}

The numerical computation for $\O=]0,1[\times]0,1[$, $f=1$, and $k=1$ gives the following outputs for various values of $m$.

\begin{figure}[htp]
\begin{center}
\includegraphics[height=6truecm]{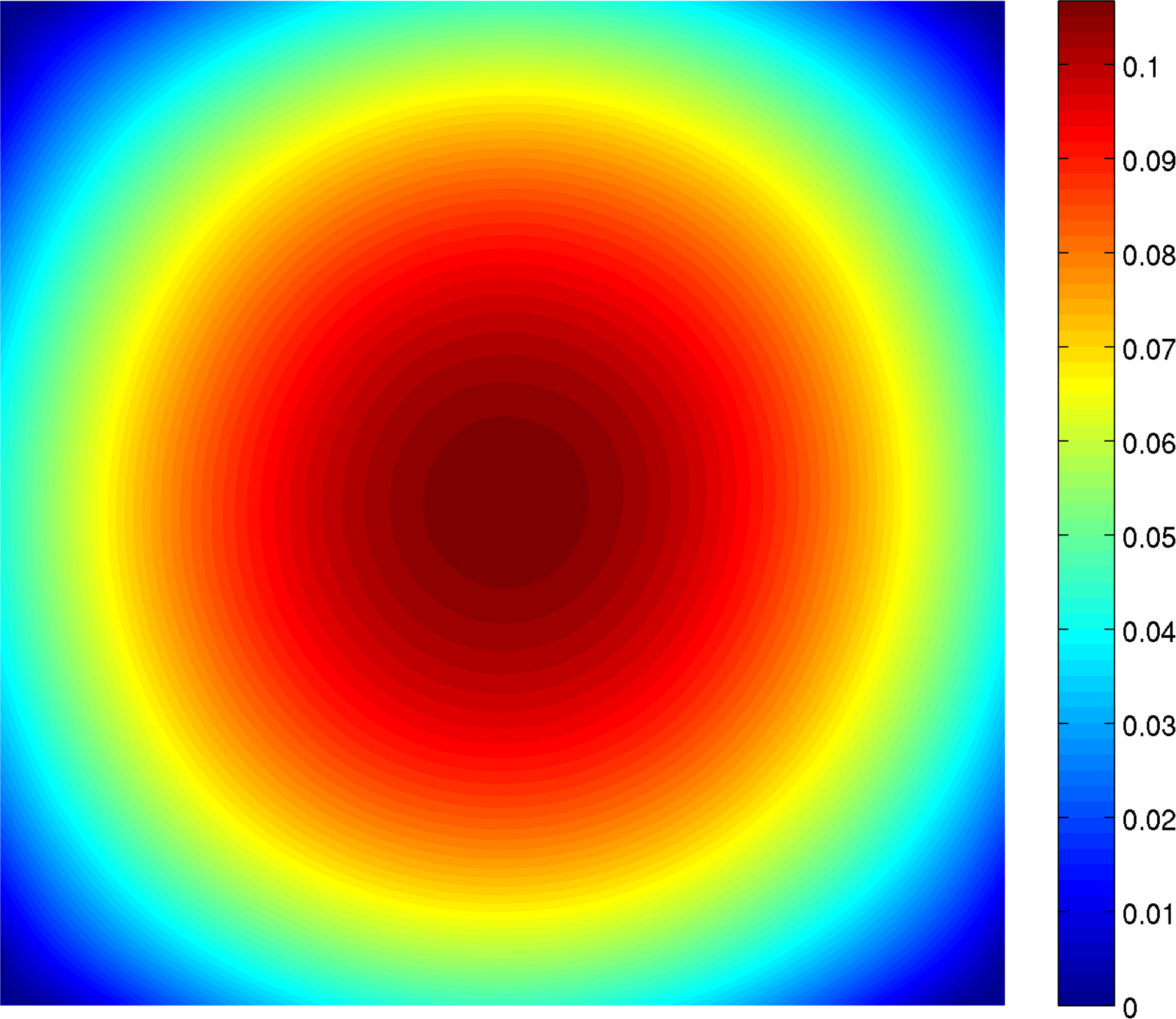}
\includegraphics[height=6truecm]{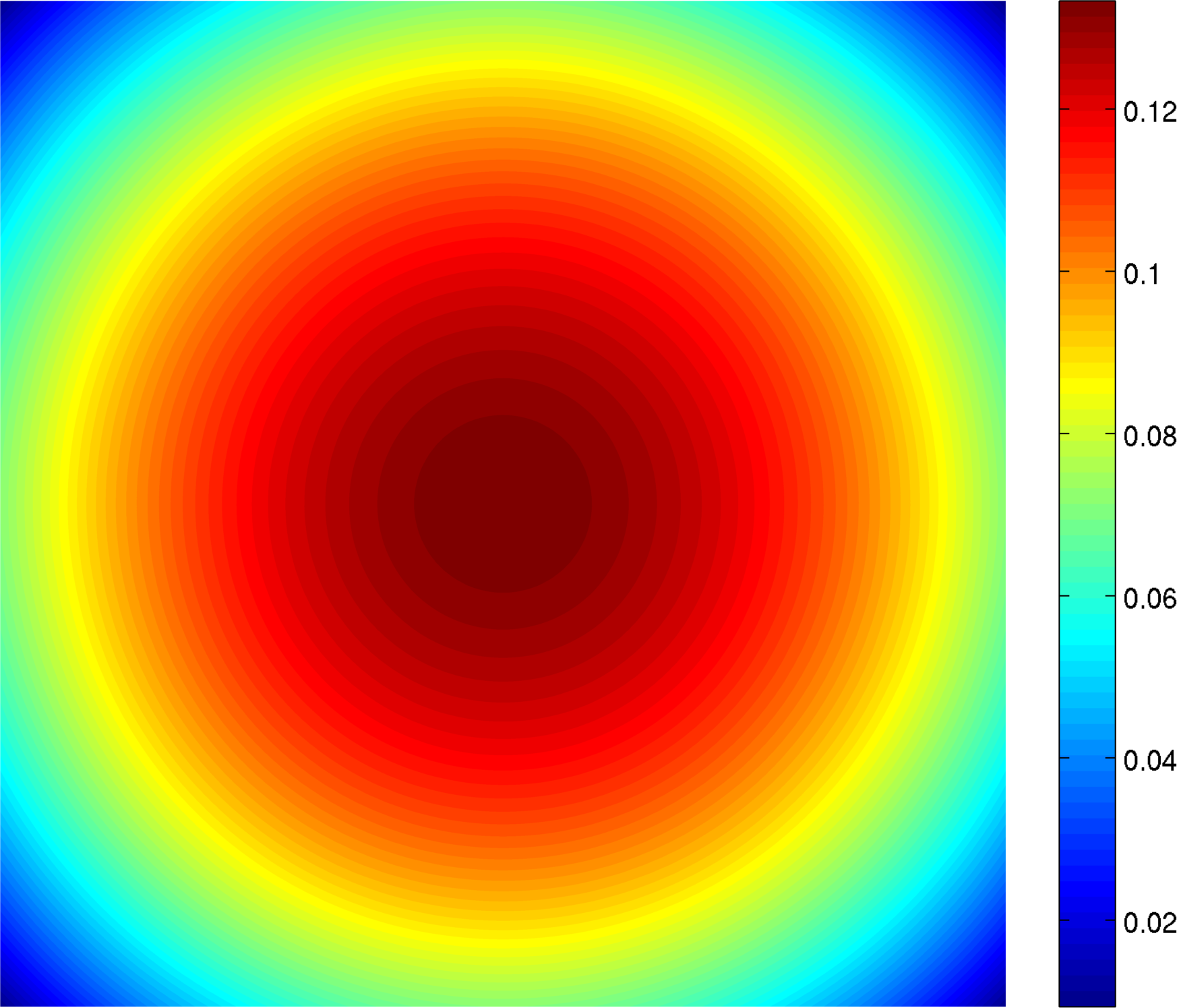}
\caption{Optimal energy insulation of $]0,1[\times]0,1[$ with $m=1$ and $m=2$.}
\label{fig4}
\end{center}
\end{figure}

It is interesting to study the behavior of $h_m$ when $m\to0$. This was analyzed in \cite{esri03}, where it is shown that the rescaled functions $h_m/m$ weakly* converge as measures to a probability measure concentrated on the set where the normal derivative $\partial u_0/\partial\nu$ reaches its minimal value, being $u_0$ the solution of the Dirichlet problem
$$-\Delta u=f,\qquad u\in H^1_0(\O).$$
For instance, when $\O$ is a square, the region of concentration of $h_m/m$ consists of the four middle points of the sides.

\section{Eigenvalue optimization}\label{seigen}

The second optimization problem \eqref{pb2} looks very similar to the first one; however, its solutions behave in a quite different way. This problem has been raised in  \cite{fr80}; in particular the question if, in the case $\O$ a ball the optimal $h$ is constant, was raised in \cite{cku99}. In that paper, the authors claimed this symmetry result; however they overlook a point in the proof, that works only if the total quantity $m$ of insulator is large enough. In \cite{bubuni} a complete proof was provided, and surprisingly, for small $m$ there is a {\it symmetry breaking} phenomenon, as specified below.

As done for problem \eqref{pb1}, problem \eqref{pb2} is a double minimization problem as well; by the same procedure of Section \ref{senergy} we find that the optimal function $h_m$ is given by
\be\label{propor}
h_m=m\frac{|\bar u|}{\int_{\partial\O}|\bar u|\,d\sigma}\;,
\ee
where $\bar u$ solves the {\it auxiliary variational problem}
\be\label{pb2aux}
\min\left\{\int_\O|\nabla u|^2\,dx+\frac{1}{km}\Big(\int_{\partial\O}|u|\,d\HH^{d-1}\Big)^2\ :\ u\in H^1(\O),\ \int_\O u^2\,dx=1\right\}.
\ee

\begin{teo}
For every $\O$ there exists a solution $h_m$ to the optimization problem \eqref{pb2}. If $\O=B_R$ there exists a threshold $m_0>0$ such that:
\begin{itemize}
\item if $m>m_0$ $\bar u$ is radial, hence $h_m$ is constant;
\item if $m<m_0$ $\bar u$ is not radial, hence $h_{opt}$ is not constant.
\end{itemize}
\end{teo}

\noindent The threshold value $m_0$ is determined as the unique $m$ such that $\lambda_m=\Lambda$, where
$$\lambda_m=\min\left\{\int_\O|\nabla u|^2\,dx+\frac{1}{km}\Big(\int_{\partial\O}|u|\,d\HH^{d-1}\Big)^2\ :\ u\in H^1(\O),\ \int_\O u^2\,dx=1\right\},$$
while $\Lambda$ is the first nonzero Neumann eigenvalue
$$\Lambda=\min\left\{\int_\O|\nabla u|^2\,dx\ :\ u\in H^1(\O),\ \int_\O u^2\,dx=1,\ \int_\O u\,dx=0\right\}.$$

When the dimension $d=1$ no symmetry breaking occurs. In fact, in this case, the first nonzero Neumann eigenvalue $\Lambda$ coincides with the first Dirichlet eigenvalue $\Lambda_0$ and so $\lambda_m<\Lambda$ for every $m$.

Below are some numerical outputs for $\O$ the unitary disc in $\R^2$ and $k=1$, for various values of $m$. The plots show the first eigenfunction, solution of the minimum problem \eqref{pb2aux}; the boundary insulator thickness is proportional to it, by \eqref{propor}.

\begin{figure}[htp]
\begin{center}
\includegraphics[scale=0.5]{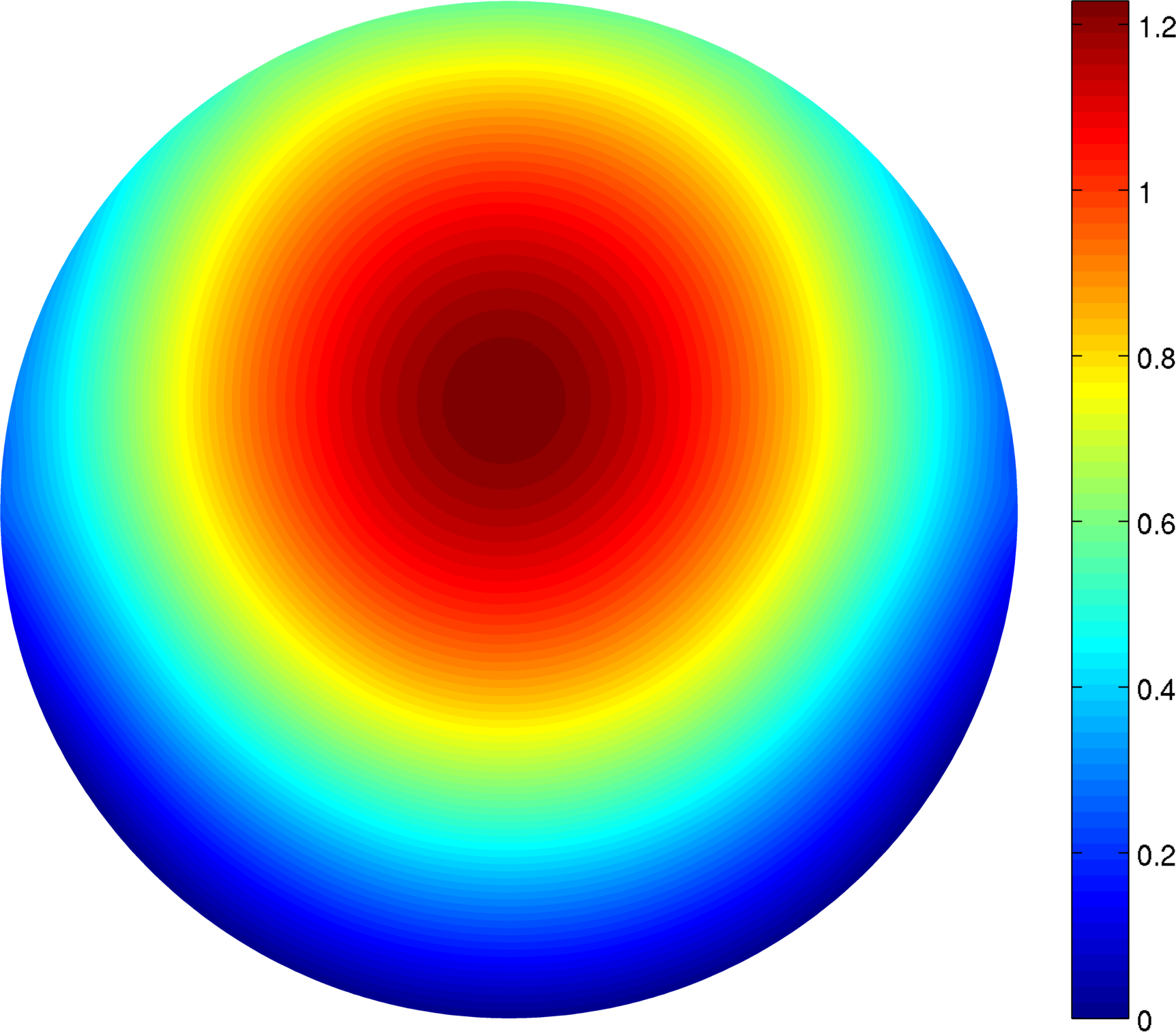}
\includegraphics[scale=0.5]{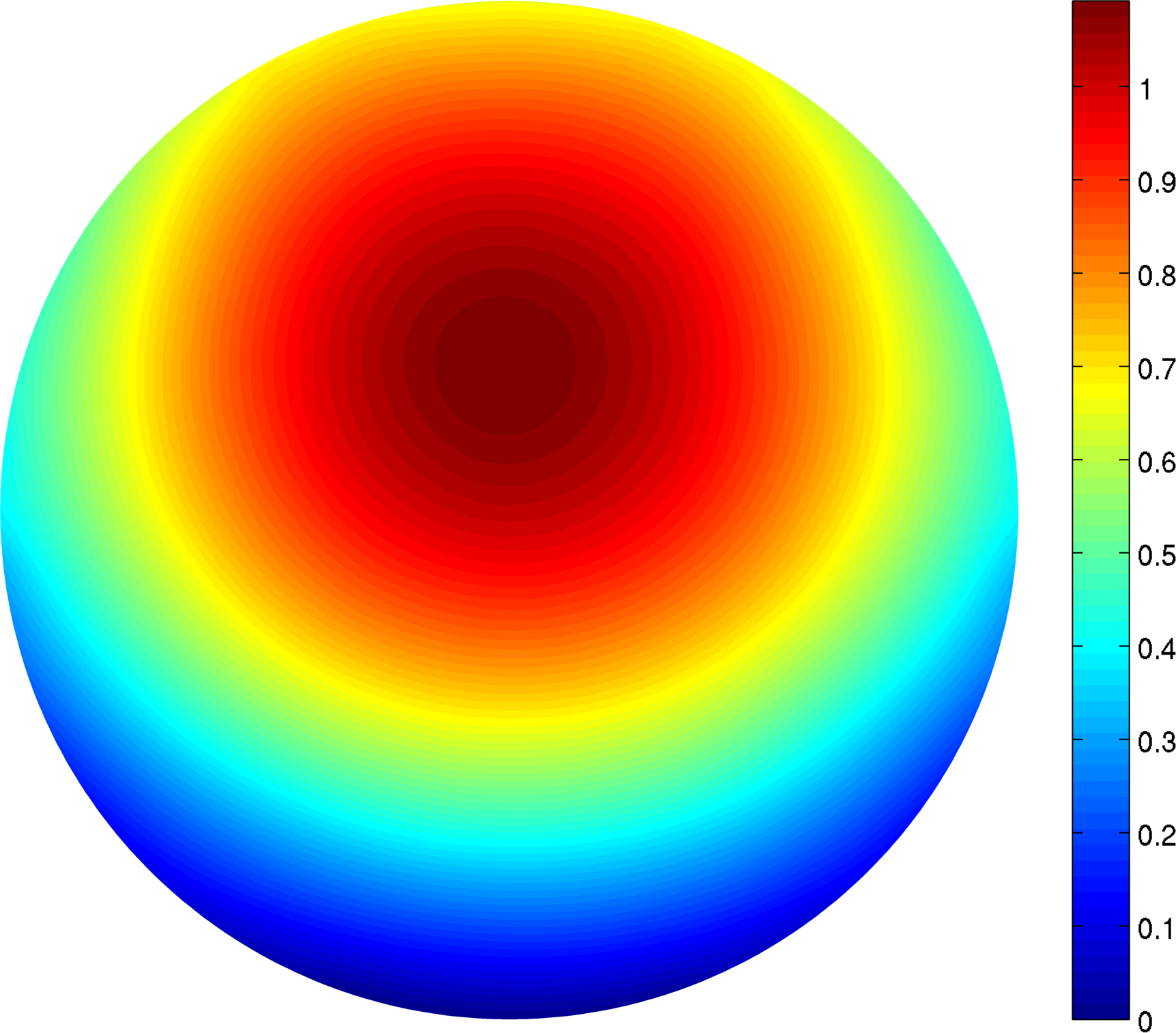}
\caption{Optimal eigenvalue insulation of the unitary disc with $k=1$, for $m=1$ and $m=2$.}
\label{fig4}
\end{center}
\end{figure}

\section{Some open problems}\label{sopen}

In the presentation above we always considered $\O$ fixed; however it would be very interesting to study the {\it shape optimization problems} related to problems 1 and 2. More precisely, denoting by $\E(\O)$ and by $\lambda(\O)$ the minimal values of problems 1 and 2
\[\begin{split}
&\E(\O)=\min\left\{\frac12\int_\O|\nabla u|^2\,dx+\frac{1}{2km}\Big(\int_{\partial\O}|u|\,d\HH^{d-1}\Big)^2-\int_\O fu\,dx\ :\ u\in H^1(\O)\right\},\\
&\lambda(\O)=\min\left\{\int_\O|\nabla u|^2\,dx+\frac{1}{km}\Big(\int_{\partial\O}|u|\,d\HH^{d-1}\Big)^2\ :\ u\in H^1(\O),\ \int_\O u^2\,dx=1\right\},
\end{split}\]
the two shape optimization problems read as
\[\begin{split}
&\min\big\{\E(\O)\ :\ |\O|=M\big\},\\
&\min\big\{\lambda(\O)\ :\ |\O|=M\big\}.
\end{split}\]

The problems above look very difficult and we do not have at the moment existence results of optimal shapes. It would be very interesting to prove (or disprove) the following facts.
\begin{itemize}
\item In the case of energy show that an optimal shape exists.

\item In the case of energy show that the optimal shape is a ball. We can actually prove that the ball is {\it stationary} with respect to smooth perturbations of the boundary.

\item In the case of eigenvalue, in dimension $d=2$, show that an optimal shape exists. If $d\ge3$ it is easy to see that, taking $\O$ as the union of many disjoint small balls, the value $\lambda(\O)$ can be made arbitrarily close to zero.

\item In the case of eigenvalue, in dimension $d=2$, the optimal shape is a ball if $m$ is large ($m>m_0$).

\item In the case of eigenvalue, in dimension $d=2$, characterize the optimal shape (if any) when $m$ is small ($m<m_0$). We can prove that if $m<m_0$ the ball cannot be optimal, actually it is even not a stationary domain with respect to smooth perturbations of the boundary.
\end{itemize}

\ack We wish to thank Beniamin Bogosel for the help provided in the numerical simulations.


\bigskip
{\small\noindent
Dorin Bucur:
Laboratoire de Math\'ematiques (LAMA),
Universit\'e de Savoie\\
Campus Scientifique,
73376 Le-Bourget-Du-Lac - FRANCE\\
{\tt dorin.bucur@univ-savoie.fr}\\
{\tt http://www.lama.univ-savoie.fr/$\sim$bucur/}

\bigskip\noindent
Giuseppe Buttazzo:
Dipartimento di Matematica,
Universit\`a di Pisa\\
Largo B. Pontecorvo 5,
56127 Pisa - ITALY\\
{\tt buttazzo@dm.unipi.it}\\
{\tt http://www.dm.unipi.it/pages/buttazzo/}

\bigskip\noindent
Carlo Nitsch:
Dipartimento di Matematica e Applicazioni,
Universit\`a di Napoli ``Federico II''\\
Via Cintia, Monte S. Angelo,
80126 Napoli - ITALY\\
{\tt carlo.nitsch@unina.it}\\
{\tt http://wpage.unina.it/c.nitsch/}

\end{document}